\documentclass[12pt]{amsart}
\usepackage{amsmath, amssymb, epsfig}

\newlength{\algorithmwidth}
\theoremstyle{plain}
\newtheorem{theorem}{Theorem}[section]
\newtheorem{proposition}[theorem]{Proposition}

\newtheorem{lemma}[theorem]{Lemma}

\theoremstyle{definition}
\newtheorem{definition}[theorem]{Definition}

\theoremstyle{remark}
\newtheorem*{remark}{Remark}

\numberwithin{equation}{section}

\DeclareMathOperator*{\supp}{supp}

\DeclareMathOperator*{\range}{range}
\DeclareMathOperator*{\Id}{Id}
\DeclareMathOperator*{\Idg}{Id_\Gamma}

\DeclareMathOperator*{\argmin}{arg min}

\def \R {\mathbb{R}}

\def \P {\mathbb{P}}

\def \e {\varepsilon}

\def \d {\delta}

\def \L {\Lambda}

\def \< {\langle}
\def \> {\rangle}
\def \^ {\widehat}

\def \supp {{\rm supp}}

\begin{document}
\title[]{Uniform Uncertainty Principle and signal recovery
  via Regularized Orthogonal Matching Pursuit}

\author{Deanna Needell
  \and Roman Vershynin}

\thanks{Partially supported by the Alfred P.~Sloan Foundation
  and by NSF DMS grants 0401032 and 0652617. \\Communicated by Emmanuel Candes.}
  
\date{July 23, 2007}

\address{Department of Mathematics,
   University of California,
   Davis, CA 95616, USA}
\email{\{dneedell,vershynin\}@math.ucdavis.edu}

\begin{abstract}
  This paper seeks to bridge the two major algorithmic approaches
  to sparse signal recovery from an incomplete set of linear measurements --
  $L_1$-mini\-mization methods and iterative methods (Matching Pursuits).
  We find a simple regularized version of Orthogonal Matching Pursuit (ROMP) 
  which has advantages of both approaches: the speed and transparency of OMP
  and the strong uniform guarantees of $L_1$-minimization.
  Our algorithm ROMP reconstructs a sparse signal in a number of iterations
  linear in the sparsity, 
  and the reconstruction is exact provided the linear measurements 
  satisfy the Uniform Uncertainty Principle.
\end{abstract}
\maketitle

\section{Introduction}

Sparse recovery problems arise in many applications ranging from medical imaging 
to error correction.
Suppose $v$ is an unknown $d$-dimensional signal with at most $n \ll d$ 
nonzero components:
$$
v \in \R^d, \qquad |\supp(v)| \le n \ll d.
$$
We call such signals $n$-sparse.
Suppose we are able to collect $N \ll d$ nonadaptive linear measurements of $v$, 
and wish to efficiently recover $v$ from these.
The measurements are given as the vector $\Phi v \in \R^N$, 
where $\Phi$ is some $N \times d$ measurement matrix.\footnote{We chose to
work with real numbers for simplicity of presentation; similar results
hold over complex numbers.}

As discussed in \cite{CDD}, exact recovery is possible with just $N=2n$. However,
recovery using only this property is not numerically feasible; the sparse recovery 
problem in general is known to be NP-hard.
Nevertheless, massive recent work in the emerging area of Compressed Sensing 
demonstrated that for several natural classes of measurement matrices $\Phi$, 
the signal $v$ can be exactly reconstructed from its measurements $\Phi v$ with
\begin{equation}                    \label{N}
  N = n \log^{O(1)} (d).
\end{equation}
In other words, the number of measurements $N \ll d$ should be almost linear in 
the sparsity $n$. 
Survey \cite{C} contains some of these results; 
the Compressed Sensing webpage \cite{CS webpage} documents progress in this area.

The two major algorithmic approaches to sparse recovery are methods
based on $L_1$-minimization and iterative methods (Matching Pursuits).
We now briefly describe these methods. Then we propose a new iterative method
that has advantages of both approaches.

\subsection{$L_1$-minimization}

This approach to sparse recovery has been advocated over decades by Donoho and his 
collaborators (see e.g. \cite{D}). 
The sparse recovery problem can be stated as the problem of finding 
the sparsest signal $v$ with the given measurements $\Phi v$:
\begin{equation}                    \tag{$L_0$}
  \min \|u\|_0 \qquad \text{subject to} \qquad \Phi u = \Phi v
\end{equation}
where $\|u\|_0 := |\supp(u)|$.
Donoho and his associates advocated the principle
that for some measurement matrices $\Phi$, the highly non-convex 
combinatorial optimization problem $(L_0)$ should be equivalent 
to its convex relaxation 
\begin{equation}                    \tag{$L_1$}
  \min \|u\|_1 \qquad \text{subject to} \qquad \Phi u = \Phi v
\end{equation}
where $\|u\|_1 = \sum_i |u_i|$ denotes the $\ell_1$-norm of the vector
$u = (u_1,\ldots,u_d)$. The convex problem $(L_1)$ can be solved using 
methods of convex and even linear programming. 

The recent progress in the emerging area of Compressed Sensing pushed forward
this program (see survey \cite{C}). A necessary and sufficient
condition of exact sparse recovery is that the map $\Phi$ be one-to-one on the set of  
$n$-sparse vectors. Cand\`es and Tao \cite{CT decoding} proved that a stronger
quantitative version of this condition guarantees the equivalence 
of the problems $(L_0)$ and $(L_1)$.

\begin{definition}[Restricted Isometry Condition] 
  A measurement matrix $\Phi$ satisfies the 
  {\em Restricted Isometry Condition} (RIC)
  with parameters $(m, \e)$ for $\e \in (0,1)$
  if we have
  $$
  (1-\e)\|v\|_2 \leq \|\Phi v\|_2 \leq (1+\e)\|v\|_2
  \qquad \text{for all $m$-sparse vectors}.
  $$
\end{definition}

The Restricted Isometry Condition states that
every set of $m$ columns of $\Phi$ forms approximately an 
orthonormal system. 
One can interpret the Restricted Isometry Condition as an abstract 
version of the Uniform Uncertainty Principle in harmonic analysis
(\cite{CT}, see also discussions in \cite{CRT} and \cite{LV}).

\begin{theorem}[Sparse recovery under RIC \cite{CT decoding}]  \label{recovery RIC}
  Assume that the measurement matrix $\Phi$ satisfies the
  Restricted Isometry Condition with parameters $(3n, 0.2)$.
  Then every $n$-sparse vector $x$ can be exactly recovered from 
  its measurements $\Phi x$ as a unique solution to the convex
  optimization problem $(L_1)$.
\end{theorem}

\medskip

In a lecture on Compressive Sampling, Cand\`es sharpened this to work for the Restricted Isometry Condition 
with parameters $(2n, \sqrt{2}-1)$. Measurement matrices that satisfy the Restricted Isometry Condition 
with number of measurements as in \eqref{N} include random Gaussian, 
Bernoulli and partial Fourier matrices. 
Section~\ref{s: ensembles} contains more detailed information.

\subsection{Orthogonal Matching Pursuit (OMP)}

An alternative approach to sparse recovery is via iterative algorithms, 
which find the support of the $n$-sparse signal $v$ 
progressively. Once $S = \supp(v)$ is found correctly, it is easy to compute the signal $v$ 
from its measurements $x = \Phi v$ as $v = (\Phi_S)^{-1} x$, 
where $\Phi_S$ denotes the measurement matrix $\Phi$ restricted to columns indexed by $S$.

A basic iterative algorithm is Orthogonal Matching Pursuit (OMP),
popularized and analyzed by Gilbert and Tropp in \cite{TG}, see \cite{T} for a more general setting.
OMP recovers the support of $v$, one index at a time, in $n$ steps. 
Under a hypothetical assumption that $\Phi$ is an isometry, i.e. 
the columns of $\Phi$ are orthonormal, 
the signal $v$ can be exactly recovered from its measurements 
$x = \Phi v$ as $v = \Phi^* x$. 

The problem is that the $N \times d$ matrix $\Phi$ 
is never an isometry in the interesting range where the number of measurements 
$N$ is smaller than the ambient dimension $d$. Even though the matrix is not an isometry, one can still
use the notion of coherence in recovery of sparse signals. In that setting, 
greedy algorithms are used with incoherent dictionaries to recover such signals, see
\cite{DET1}, \cite{DET2}, \cite{GMS}. In our setting, for random matrices
one expects the columns to be approximately orthogonal, 
and the {\em observation vector} 
$u = \Phi^* x$ to be a good approximation to the original signal $v$. 

The biggest coordinate of the observation vector $u$ in magnitude should thus 
be a nonzero coordinate of the signal $v$. We thus find one point of the support of $v$. 
Then OMP can be described as follows. First, we initialize the residual $r = x$.
At each iteration, we compute the observation vector $u = \Phi^* r$. 
Denoting by $I$ the coordinates selected so far, we solve
a least squares problem and update the residual
$$
y = \argmin_{z \in \R^I} \|x - \Phi z\|_2; \qquad r = x - \Phi y,
$$
to remove any contribution of the coordinates in $I$.
OMP then iterates this procedure $n$ times, and outputs a set $I$
of size $n$, which should equal the support of the signal $v$. 

Tropp and Gilbert \cite{TG} analyzed the performance of OMP for Gaussian 
measurement matrices $\Phi$; a similar result holds for general subgaussian matrices.
They proved that, for every fixed $n$-sparse $d$-dimensional signal $v$, and 
an $N \times d$ random Gaussian measurement matrix $\Phi$, OMP recovers 
(the support of) $v$ from the measurements $x = \Phi v$ correctly with high 
probability, provided the number of measurements is $N \sim n \log d$. 

\subsection{Advantages and challenges of both approaches}

The $L_1$-mini\-mization method has {\em strongest known guarantees} of sparse recovery. 
Once the measurement matrix $\Phi$ satisfies the Restricted Isometry Condition, 
this method works correctly for all sparse signals $v$. 
No iterative methods have been known to feature 
such uniform guarantees, with the exception of Chaining Pursuit \cite{GSTV} 
and the HHS Algorithm \cite{HHS}
which however only work with specifically designed structured
measurement matrices. 

The Restricted Isometry Condition is a natural abstract deterministic property 
of a matrix. Although establishing this property is often nontrivial, 
this task is {\em decoupled from the analysis} of the recovery algorithm.

$L_1$-minimization is based on linear programming, which has its advantages
and disadvantages. One thinks of linear programming as a black box, 
and any development of fast solvers will reduce the {\em running time} of 
the sparse recovery method. On the other hand, it is not very clear what this running time is, 
as there is no strongly polynomial time algorithm 
in linear programming yet. All known solvers take time polynomial not only in 
the dimension of the program $d$, but also on certain condition numbers of the program.
While for some classes of random matrices the expected running time of linear 
programming solvers can be bounded (see the discussion in \cite{ST} and 
subsequent work in \cite{V}), estimating
condition numbers is hard for specific matrices. For example, there is no result
yet showing that the Restricted Isometry Condition implies that the condition numbers
of the corresponding linear program is polynomial in $d$.

Orthogonal Matching Pursuit is quite {\em fast}, both theoretically and experimentally.
It makes $n$ iterations, where each 
iteration amounts to a multiplication by a $d \times N$ matrix $\Phi^*$
(computing the observation vector $u$), and solving a least squares problem 
in dimensions at most $N \times n$ (with matrix $\Phi_I$). This yields strongly 
polynomial running time. In practice, OMP is observed to perform faster and 
is easier to implement than $L_1$-minimization \cite{TG}.
For more details, see \cite{TG}.

Orthogonal Matching Pursuit is quite {\em transparent}: at each iteration, it
selects a new coordinate from the support of the signal $v$ in a very specific and 
natural way. In contrast, the known $L_1$-minimization solvers, such as the 
simplex method and interior point methods, compute a path toward the solution. 
However, the geometry of $L_1$ is clear, whereas the analysis of greedy algorithms
can be difficult simply because they are iterative.

On the other hand, Orthogonal Matching Pursuit has {\em weaker guarantees} of 
exact recovery. Unlike $L_1$-minimization, the guarantees of OMP are non-uniform:
for each {\em fixed} sparse signal $v$ and not for {\em all} signals, 
the algorithm performs correctly with high probability. 
Rauhut has shown that uniform guarantees for OMP are impossible for natural 
random measurement matrices \cite{R}.

Moreover, OMP's condition on measurement matrices given in \cite{TG} 
is {\em more restrictive} than the Restricted Isometry Condition. 
In particular, it is not known whether OMP succeeds in the important 
class of partial Fourier measurement matrices.

These open problems about OMP, first stated in \cite{TG} and often reverberated 
in the Compressed Sensing community, motivated the present paper.
We essentially settle them in positive by the following modification of 
Orthogonal Matching Pursuit.

\subsection{Regularized OMP}

This new algorithm for sparse recovery will perform correctly for all 
measurement matrices $\Phi$ satisfying the Restricted Isometry Condition,
and for all sparse signals.

When we are trying to recover the signal $v$ from its measurements $x = \Phi v$,
we can use the observation vector $u = \Phi^* x$ as a good {\em local 
approximation} to the signal $v$. Namely, the observation vector $u$
encodes correlations of the measurement vector $x$ with the columns of $\Phi$. Note that $\Phi$ is a 
dictionary, and so since the signal $v$ is sparse, $x$ has a sparse representation with
respect to the dictionary. By the Restricted 
Isometry Condition, every $n$ columns form approximately an orthonormal system.
Therefore, every $n$ coordinates of the observation vector $u$ look like correlations of the measurement vector $x$ 
with the orthonormal basis and therefore are close in the Euclidean norm to the corresponding $n$ coefficients of $v$. 
This is documented in Proposition~\ref{P:cons} below.

The local approximation property suggests to make use of the $n$ biggest coordinates of 
the observation vector $u$, rather than one biggest coordinate as OMP did.
We thus force the selected coordinates to be more regular (ie. closer to uniform) by selecting only the coordinates with comparable sizes. To this end, a new {\em regularization} step will be needed to ensure
that each of these coordinates gets an even share of information. This leads to the following algorithm for sparse recovery:

\bigskip

\textsc{Regularized Orthogonal Matching Pursuit (ROMP)}

\nopagebreak

\fbox{\parbox{\algorithmwidth}{
  \textsc{Input:} Measurement vector $x \in \R^N$ and sparsity level $n$
  
  \textsc{Output:} Index set $I \subset \{1,\ldots,d\}$

  \begin{description}
    \item[Initialize] Let the index set $I = \emptyset$ and the residual $r = x$.\\
      Repeat the following steps until $r = 0$:
    \item[Identify] Choose a set $J$ of the $n$ biggest coordinates in magnitude 
      of the observation vector $u = \Phi^*r$, or all of its nonzero coordinates, 
      whichever set is smaller.
    \item[Regularize] Among all subsets $J_0 \subset J$ with comparable coordinates:
      $$
      |u(i)| \leq 2|u(j)| \quad \text{for all } i,j \in J_0,
      $$
      choose $J_0$ with the maximal energy $\|u|_{J_0}\|_2$.
    \item[Update] Add the set $J_0$ to the index set: $I \leftarrow I \cup J_0$, 
      and update the residual:
      $$
      y = \argmin_{z \in \R^I} \|x - \Phi z\|_2; \qquad r = x - \Phi y.
      $$
  \end{description}
 }}

    \bigskip
  
  \begin{remark} The identification and regularization steps of ROMP can be performed
  efficiently. In particular, the regularization step does \textit{not} imply combinatorial 
  complexity, but actually can be done in linear time. 
  The running time of ROMP is thus comparable to that of OMP in theory, 
  and is often better than OMP in practice. We discuss the runtime in detail in Section~\ref{s: implementation}.
  \end{remark}

The main theorem of this paper states that ROMP yields exact sparse recovery
provided that the measurement matrix satisfies the Restricted Isometry Condition.

\begin{theorem}[Exact sparse recovery via ROMP]\label{T:main}
  Assume a measurement matrix $\Phi$ satisfies the Restricted Isometry Condition 
  with parameters $(2n, \e)$ for $\e = 0.03 / \sqrt{\log n}$. 
  Let $v$ be an $n$-sparse vector in $\R^d$ with measurements $x = \Phi v$. 
  Then ROMP in at most $n$ iterations outputs a set $I$ such that 
  $$
  \supp(v) \subset I \quad \text{and} \quad |I| \leq 2n.
  $$
\end{theorem}

\medskip

This theorem is proved in Section \ref{s: proof}.

\medskip

\begin{remarks}
  {\bf 1. } Theorem~\ref{T:main} guarantees {\em exact sparse recovery}.
  Indeed, it is easy to compute the signal $v$ from its measurements $x = \Phi v$ 
  and the set $I$ given by ROMP as $v = (\Phi_I)^{-1} x$,
  where $\Phi_I$ denotes the measurement matrix $\Phi$ restricted to columns
  indexed by $I$.

  {\bf 2.} Theorem~\ref{T:main} gives {\em uniform guarantees} of sparse recovery.
  Indeed, once the measurement matrix satisfies a deterministic condition (RIC),
  then our algorithm ROMP correctly recovers {\em every} sparse vector from its measurements. 
  Uniform guarantees have been shown to be impossible for OMP \cite{R}, and it has been 
  an open problem to find a version of OMP with uniform guarantees (see \cite{TG}). 
  Theorem~\ref{T:main} says that ROMP essentially settles this problem.

  {\bf 3. } The logarithmic factor in $\e$ may be an artifact of the proof. 
  At this moment, we do not know how to remove it.
    
  {\bf 4. } Measurement matrices known to satisfy the Restricted Isometry 
  Condition include random {\em Gaussian, Bernoulli and partial Fourier matrices},
  with number of measurements $N$ almost linear in the sparsity $n$, i.e. 
  as in \eqref{N}. Section~\ref{s: ensembles} contains detailed information.
  It has been unknown whether OMP gives sparse recovery for partial Fourier
  measurements (even with non-uniform guarantees). ROMP gives sparse recovery 
  for these measurements, and even with uniform guarantees. 
  
\end{remarks}

\medskip

The rest of the paper is organized as follows. 
In Section~\ref{s: ensembles} we describe known classes of measurement matrices
satisfying the Restricted Isometry Condition.
In Section~\ref{s: proof} we give the proof of Theorem~\ref{T:main}.
In Section~\ref{s: implementation} we discuss implementation, running time, and
empirical performance of ROMP.

\subsection*{Acknowledgment}
We would like to thank the referees for a thorough reading of the manuscript and making
useful suggestions which greatly improved the paper.

\section{Measurement matrices satisfying the Restricted Isometry Condition} 
\label{s: ensembles}

The only known measurement matrices known to satisfy the Restricted Isometry
Condition with number of measurements as in \eqref{N} are certain classes
of random matrices. The problem of deterministic constructions is still open. 
The known classes include: subgaussian random matrices (in particular, Gaussian
and Bernoulli), and random partial bounded orthogonal matrices (in particular, partial 
Fourier matrices).

Throughout the paper, $C, c, C_1, C_2, c_1, c_2, \ldots$ denote positive 
absolute constants unless otherwise specified.

A {\em subgaussian random matrix} $\Phi$ is a matrix whose entries are i.i.d.
subgaussian random variables with variance $1$.
A random variable $X$ is subgaussian if its tail distribution is dominated by
that of the standard Gaussian random variable: there are constants $C_1, c_1 > 0$ such
that $\P(|X| > t) \le C_1 \exp(-c_1 t^2)$ for all $t > 0$.
Examples of subgaussian random variables are: standard Gaussian,
Bernoulli (uniform $\pm 1$), and any bounded random variables.

A {\em partial bounded orthogonal matrix} $\Phi$ is formed by $N$ randomly uniformly
chosen rows of an orthogonal $d \times d$ matrix $\Psi$, whose entries
are bounded by $C_2/\sqrt{d}$, for some constant $C_2$. An example of $\Psi$
is the discrete Fourier transform matrix.
Taking measurements $\Phi v$ with a partial Fourier matrix
thus amounts to observing $N$ random frequencies of the signal $v$.

The following theorem documents known results on the Restricted Isometry Condition
for these classes of random matrices. 

\begin{theorem}[Measurement matrices satisfying RIC]    \label{RIC matrices}
  Consider an $N \times d$ measurement matrix $\Phi$, and let $n \ge 1$,
  $\e \in (0,1/2)$, and $\d \in (0,1)$.

  1. If $\Phi$ is a {\em subgaussian matrix}, then with probability $1-\d$ the matrix
  $\frac{1}{\sqrt{N}} \Phi$ satisfies the Restricted
  Isometry Condition with parameters $(n,\e)$ provided that
  $$
  N \ge \frac{Cn}{\e^2} \log \Big( \frac{d}{\e^2 n} \Big).
  $$
  2. If $\Phi$ is a {\em partial bounded orthogonal matrix},
  then with probability $1-\d$ the matrix
  $\sqrt{\frac{d}{N}} \, \Phi$ satisfies the Restricted Isometry Condition
  with parameters $(n,\e)$ provided that
  $$
  N \ge C \Big( \frac{n \log d}{\e^2} \Big)
  \log \Big( \frac{n \log d}{\e^2} \Big)
  \log^2 d.
  $$
  In both cases, the constant $C$ depends only on the confidence level 
  $\d$ and the constants $C_1, c_1, C_2$ from the definition of the corresponding 
  classes of matrices.
\end{theorem}

\begin{remarks}
  {\bf 1. } The first part of this theorem is proved in \cite{MPT}.
    The second part is from \cite{RV}; a similar estimate with somewhat worse exponents
    in the logarithms was proved in \cite{CT}. See these results for the exact
    dependence of $C$ on the confidence level $\d$ (although usually $\d$ would be chosen
    to be some small constant itself.) 
    
  {\bf 2. } In Theorem~\ref{T:main}, we needed to use RIC for $\e = c_1 / \sqrt{\log n}$. 
    An immediate consequence of Theorem~\ref{RIC matrices} is that 
    subgaussian matrices satisfy such RIC for the number of measurements
    $$
    N \sim n \log^2 d
    $$
    and partial bounded orthogonal matrices for
    $$
    N \sim n \log^5 d.
    $$
    These numbers of measurements guarantee exact sparse recovery using ROMP.
\end{remarks}

\section{Proof of Theorem~\ref{T:main}}
\label{s: proof}

We shall prove a stronger version of Theorem~\ref{T:main}, which states that 
{\em at every iteration} of ROMP, at least $50\%$ of the newly selected coordinates
are from the support of the signal $v$. 

\begin{theorem}[Iteration Invariant of ROMP] \label{T:it}
  Assume $\Phi$ satisfies the Restricted Isometry Condition 
  with parameters $(2n, \e)$ for $\e = 0.03 / \sqrt{\log n}$. 
  Let $v \ne 0$ be an $n$-sparse vector with measurements $x = \Phi v$. 
  Then at any iteration of ROMP, after the regularization step, we have
  $J_0 \ne \emptyset$, $J_0 \cap I = \emptyset$ and 
  \begin{equation}                  \label{J support}
    |J_0 \cap \supp(v)| \geq \frac{1}{2}|J_0|.
  \end{equation}
  In other words, at least $50\%$ of the coordinates in the newly selected set $J_0$ 
  belong to the support of $v$.
\end{theorem}

In particular, at every iteration ROMP finds at least one new coordinate 
in the support of the signal $v$. Coordinates outside the support can also be found, 
but \eqref{J support} guarantees that the number of such ``false'' coordinates 
is always smaller than those in the support. 
This clearly implies Theorem~\ref{T:main}.

\medskip

Before proving Theorem~\ref{T:it} we explain how the Restricted Isometry Condition 
will be used in our argument. RIC is necessarily a local principle, which concerns 
not the measurement matrix $\Phi$ as a whole, but its submatrices of $n$ columns. 
All such submatrices $\Phi_I$, $I \subset \{1,\ldots,d\}$, $|I| \le n$ are almost isometries.
Therefore, for every $n$-sparse signal $v$, the observation vector 
$u = \Phi^*\Phi v$ approximates $v$ locally, when restricted to a set of 
cardinality $n$. The following proposition formalizes these local properties
of $\Phi$ on which our argument is based.

\begin{proposition}[Consequences of Restricted Isometry Condition]\label{P:cons} 
  Assume a measurement matrix $\Phi$ satisfies the Restricted Isometry Condition 
  with parameters $(2n, \e)$. Then the following holds.
  \begin{enumerate}
    \item {\em (Local approximation)} 
      For every $n$-sparse vector $v \in \R^d$ 
      and every set $I \subset \{1, \ldots, d\}$, $|I| \le n$, 
      the observation vector $u = \Phi^* \Phi x$ satisfies 
      $$
      \|u|_I - v|_I\|_2 \le 2.03 \e \|v\|_2.
      $$
    \item {\em (Spectral norm)} 
      For any vector $z \in \R^N$
      and every set $I \subset \{1, \ldots, d\}$, $|I| \le 2n$, we have 
      $$
      \|(\Phi^*z)|_I\|_2 \leq (1+\e)\|z\|_2.
      $$
    \item {\em (Almost orthogonality of columns)} 
      Consider two disjoint sets $I,J \subset \{1, \ldots, d\}$, $|I \cup J| \le 2n$.
      Let $P_I, P_J$ denote the orthogonal projections in $\R^N$
      onto $\range(\Phi_I)$ and $\range(\Phi_J)$, respectively. Then 
      $$
      \|P_I P_J\|_{2\rightarrow 2} \leq 2.2 \e.
      $$
  \end{enumerate}
\end{proposition}

\begin{proof}

{\sc Part 1.}
Let $\Gamma = I\cup\supp(v)$, so that $|\Gamma| \leq 2n$. 
Let $\Idg$ denote the identity operator on $\R^\Gamma$. 
By the Restricted Isometry Condition,
$$
\|\Phi_\Gamma^*\Phi_\Gamma - \Idg\|_{2\rightarrow 2} 
= \sup_{y \in \R^\Gamma, \, \|y\|_2 = 1} \big| \|\Phi_\Gamma y\|_2^2 - \|y\|_2^2 \big| 
\le (1+\e)^2 - 1
\le 2.03\e.
$$
Since $\supp(v) \subset \Gamma$, we have
$$
\|u|_\Gamma - v|_\Gamma\|_2 
= \|\Phi_\Gamma^*\Phi_\Gamma v - \Id_\Gamma v\|_2 
\le 2.03\e \|v\|_2.
$$   
The conclusion of Part~1 follows since $I\subset \Gamma$.

\medskip

{\sc Part 2.} Denote by $Q_I$ the orthogonal projection in $\R^d$ 
onto $\R^I$. Since $|I| \le 2n$, the Restricted Isometry Condition yields
$$    
\|Q_I\Phi^*\|_{2\rightarrow 2} = \|\Phi Q_I\|_{2 \rightarrow 2} \leq 1 + \e.
$$
This yields the inequality in Part 2. 

\medskip

{\sc Part 3.} 
The desired inequality is equivalent to:
$$
\frac{|\langle x, y\rangle|}{\|x\|_2\|y\|_2} \leq 2.2\e \qquad \text{for all } x \in \range(\Phi_I), \; y \in \range(\Phi_J).
$$
Let $K = I \cup J$ so that $|K| \leq 2n$. For any $x \in \range(\Phi_I), y \in \range(\Phi_J)$, there are $a, b$ so that
$$
x = \Phi_K a, \; y = \Phi_K b, \qquad a \in \R^I, \; b \in \R^J.
$$
By the Restricted Isometry Condition,
$$
\|x\|_2 \geq (1-\e)\|a\|_2, \; \|y\|_2 \geq (1-\e)\|b\|_2.
$$
By the proof of Part 2 above and since $\langle a b \rangle = 0$, we have
$$
|\langle x, y\rangle| = |\langle (\Phi_K^* \Phi_K - \Idg)a, b \rangle| \leq 2.03\e\|a\|_2\|b\|_2.
$$
This yields  
$$
\frac{|\langle x, y\rangle|}{\|x\|_2\|y\|_2} \leq \frac{2.03\e}{(1-\e)^2} \leq 2.2\e,
$$
which completes the proof.
\end{proof}

\bigskip

We are now ready to prove Theorem \ref{T:it}.

\medskip

The proof is by induction on the iteration of ROMP.
The induction claim is that for all previous iterations, the set of newly chosen 
indices $J_0$ is nonempty, disjoint from the set of previously chosen indices $I$, 
and \eqref{J support} holds.
 
Let $I$ be the set of previously chosen indices at the start of a given iteration.
The induction claim easily implies that 
\begin{equation}					\label{suppv I}
  |\supp(v) \cup I| \le 2n.    
\end{equation}
Let $J_0$, $J$, be the sets found by ROMP in the current iteration. 
By the definition of the set $J_0$, it is nonempty.

Let $r \ne 0$ be the residual at the start of this iteration. 
We shall approximate $r$ by a vector in $\range(\Phi_{\supp(v) \setminus I})$.
That is, we want to approximately realize the residual $r$  
as measurements of some signal which lives on the still unfound
coordinates of the the support of $v$. 
To that end, we consider the subspace
$$
H := \range (\Phi_{\supp(v) \cup I})
$$
and its complementary subspaces
$$
F := \range (\Phi_I), \quad 
E_0 := \range (\Phi_{\supp(v) \setminus I}).
$$
The Restricted Isometry Condition in the form of Part~3 of Proposition~\ref{P:cons} 
ensures that $F$ and $E_0$ are almost orthogonal. Thus $E_0$ is close to 
the orthogonal complement of $F$ in $H$,
$$
E := F^{\perp}\cap H.
$$
\begin{figure}[ht] \label{fig:cap}
  \includegraphics[scale=0.6]{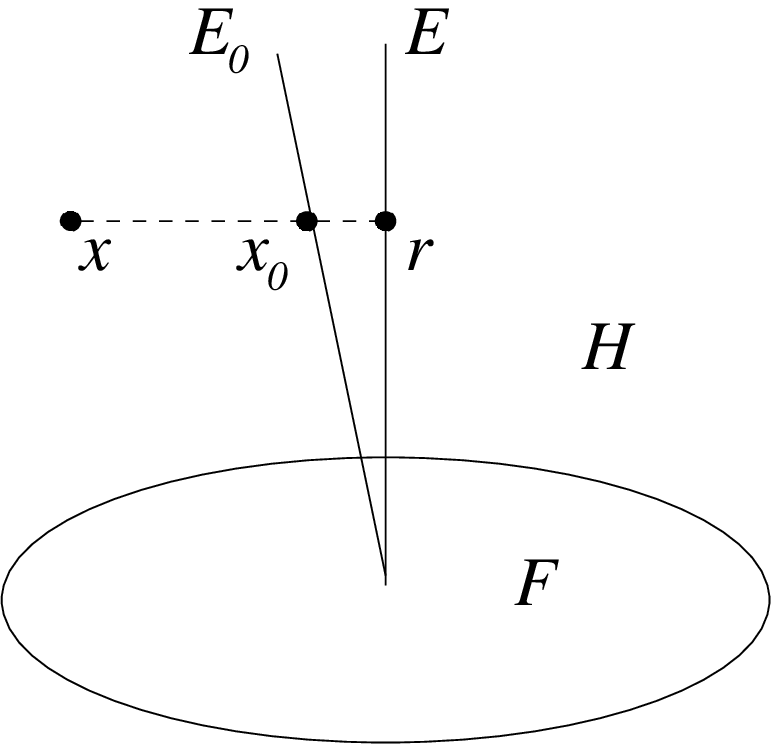}
\end{figure}

We will also consider the signal we seek to identify at the current iteration, 
its measurements, and its observation vector:
\begin{equation}            \label{v0 x0}
  v_0 := v|_{\supp(v) \setminus I}, \quad 
  x_0 := \Phi v_0 \in E_0, \quad u_0 := \Phi^*x_0.
\end{equation}

Lemma~\ref{L:uj} will show that $\|(u-u_0)|_T\|_2$ for any small enough subset $T$ is small, and Lemma~\ref{C:uj0} will show that $\|u|_{J_0}\|_2$ is not too small. First, we show that the residual $r$ has a simple description:

\begin{lemma}[Residual]     \label{residual}
  Here and thereafter, let $P_L$ denote the orthogonal projection in $\R^N$ 
  onto a linear subspace $L$. Then
  $$
  r = P_E x.
  $$
\end{lemma}

\begin{proof}
By definition of the residual in the algorithm,
$r = P_{F^\perp} x$. Since $x \in H$, we conclude from the orthogonal 
decomposition $H = F + E$ that $x = P_F x + P_E x$. Thus
$r = x - P_F x = P_E x$.
\end{proof}

To guarantee a correct identification of $v_0$, we first state
two approximation lemmas that reflect in two different ways the fact 
that subspaces $E_0$ and $E$ are close to each other.
This will allow us to carry over information from $E_0$ to $E$.

\begin{lemma}[Approximation of the residual]\label{C:proj}
  We have
  $$
  \|x_0 - r\|_2 \leq 2.2 \e \|x_0\|_2.
  $$
\end{lemma}

\begin{proof}
By definition of $F$, we have 
$x - x_0 = \Phi(v - v_0) \in F$. 
Therefore, by Lemma~\ref{residual},
$r = P_{E}x = P_{E}x_0$, and so
$$
x_0 - r = x_0 - P_Ex_0 = P_Fx_0 = P_FP_{E_0}x_0.
$$
Now we use Part 3 of Proposition~\ref{P:cons} for the sets $I$ and $\supp(v) \setminus I$
whose union has cardinality at most $2n$ by \eqref{suppv I}. It follows that
$\|P_FP_{E_0}x_0\|_2 \le 2.2 \e \|x_0\|_2$ as desired.
\end{proof}

\begin{lemma}[Approximation of the observation]\label{L:uj}
  Consider the observation vectors 
  $u_0 = \Phi^*x_0$ and $u = \Phi^*r$. Then for any set $T \subset \{1, \ldots, d\}$ with $|T| \le 2n$, we have
  $$
  \|(u_0 - u)|_T\|_2 \leq 2.4 \e \|v_0\|_2.
  $$
\end{lemma}

\begin{proof}
Since $x_0 = \Phi v_0$, we have by Lemma~\ref{C:proj} 
and the Restricted Isometry Condition that
$$
\|x_0 - r\|_2 
\le 2.2 \e \|\Phi v_0\|_2 
\le 2.2 \e (1+\e) \|v_0\|_2 
\le 2.3 \e \|v_0\|_2.
$$
To complete the proof, it remains to apply Part 2 of Proposition~\ref{P:cons},
which yields 
$\|(u_0 - u)|_T\|_2 \le (1 + \e)\|x_0 - r\|_2$.
\end{proof}

We next show that the energy (norm) of $u$ when restricted to $J$, and furthermore to 
$J_0$, is not too small. By the approximation lemmas, this will yield that ROMP 
selects at least a fixed percentage of energy of the still unidentified part of the signal. 
By the regularization step of ROMP, since all selected coefficients have comparable
magnitudes, we will conclude that not only a portion of energy
but also of the {\em support} is selected correctly. This 
will be the desired conclusion.

\begin{lemma}[Localizing the energy]\label{C:uj}
  We have $\|u|_J\|_2 \ge 0.8 \|v_0\|_2$.
\end{lemma}

\begin{proof}
Let $S$ = $\supp(v) \setminus I$.
Since $|S| \leq n$, the maximality property of $J$ in the algorithm
implies that 
$$
\|u_0|_J\|_2 \geq \|u_0|_S\|_2.
$$
Furthermore, since $v_0|_S = v_0$, by Part 1 of Proposition~\ref{P:cons} we have
$$
\|u_0|_S\|_2 \geq (1 - 2.03\e)\|v_0\|_2.
$$
Putting these two inequalities together and using
Lemma~\ref{L:uj}, we conclude that
$$
\|u|_J\|_2 \ge (1 - 2.03\e)\|v_0\|_2 - 2.4\e\|v_0\|_2 \ge 0.8 \|v_0\|_2.
$$
This proves the lemma.
\end{proof}

We next bound the norm of $u$ restricted to the smaller set $J_0$. 
We do this by first noticing a general property of regularization:

\begin{lemma}[Regularization]\label{L:reg}
  Let $y$ be any vector in $\R^m$, $m > 1$.
  Then there exists a subset $A \subset \{1, \ldots, m\}$ 
  with comparable coordinates:
  \begin{equation}\label{E:comp}
    |y(i)| \le 2|y(j)| \quad \text{for all $i, j \in A$,}
  \end{equation}
  and with big energy:
  \begin{equation}\label{E:big} 
    \|y|_A\|_2 \ge \frac{1}{2.5\sqrt{\log m}}\|y\|_2.
  \end{equation}
\end{lemma}

\begin{proof}
We will construct at most $O(\log m)$ subsets $A_k$ with comparable coordinates 
as in (\ref{E:comp}), and such that at least one of these sets will have 
large energy as in (\ref{E:big}). 

Let $y = (y_1, \ldots, y_m)$, and consider a partition of $\{1, \ldots, m\}$ 
using sets with comparable coordinates:
$$
A_k := \{i : 2^{-k}\|y\|_2 < |y_i| \leq 2^{-k+1}\|y\|_2 \}, \qquad k=1, 2, \ldots 
$$
Let $k_0 = \left\lceil \log m \right\rceil + 1$, so that 
$|y_i| \leq \frac{1}{m}\|y\|_2$ for all $i \in A_k$, $k > k_0$.
Then the set $U = \bigcup_{k \le k_0} A_k$ contains most of the energy of $y$:
$$
\|y|_{U^c}\|_2 \le \big(m(\frac{1}{m}\|y\|_2)^2\big)^{1/2} = \frac{1}{\sqrt{m}}\|y\|_2
\le \frac{1}{\sqrt{2}}\|y\|_2.
$$
Thus
$$
\big( \sum_{k\leq k_0}\|y|_{A_k}\|_2^2 \big)^{1/2} 
= \|y|_U\|_2 = \big(\|y\|_2^2 - \|y|_{U^c}\|_2^2\big)^{1/2}
\ge \frac{1}{\sqrt{2}} \|y\|_2.
$$
Therefore there exists $k \leq k_0$ such that
$$
\|y|_{A_k}\|_2 \ge \frac{1}{\sqrt{2k_0}} \|y\|_2 
\ge \frac{1}{2.5\sqrt{\log m}}\|y\|_2,
$$
which completes the proof.
\end{proof}

In our context, Lemma~\ref{L:reg} applied to the vector $u|_J$ along with 
Lemma~\ref{C:uj} directly implies:

\begin{lemma}[Regularizing the energy]\label{C:uj0} 
  We have
  $$
  \|u|_{J_0}\|_2 \ge \frac{0.32}{\sqrt{\log n}}\|v_0\|_2.
  $$
\end{lemma}

\medskip

We now finish the proof of Theorem~\ref{T:it}. 

To show the first claim, that $J_0$ is nonempty, we note that 
$v_0 \ne 0$. Indeed, otherwise by \eqref{v0 x0} we have $I \subset \supp(v)$, 
so by the definition of the residual in the algorithm, we would have $r = 0$
at the start of the current iteration, which is a contradiction. 
Then $J_0 \ne \emptyset$ by Lemma~\ref{C:uj0}.

The second claim, that $J_0 \cap I = \emptyset$, is also simple. 
Indeed, recall that by the definition of the algorithm, 
$r = P_{F^\perp} \in F^\perp = (\range(\Phi_I))^\perp$. 
It follows that the observation vector $u = \Phi^* r$  
satisfies $u|_I = 0$. Since by its definition the set $J$ contains only 
nonzero coordinates of $u$ we have $J \cap I = \emptyset$.
Since $J_0 \subset J$, the second claim $J_0 \cap I = \emptyset$ follows. 

The nontrivial part of the theorem is its last claim, inequality \eqref{J support}.
Suppose it fails. Namely, suppose that 
$|J_0 \cap \supp(v)| < \frac{1}{2}|J_0|$, 
and thus
$$
|J_0 \backslash \supp(v)| > \frac{1}{2}|J_0|.
$$
Set $\Lambda = J_0\backslash\supp(v)$. 
By the comparability property of the coordinates in $J_0$ 
and since $|\Lambda| > \frac{1}{2}|J_0|$, there is a fraction of energy 
in $\Lambda$:
\begin{equation}\label{E:ubig} 
  \|u|_{\Lambda}\|_2 > \frac{1}{\sqrt{5}}\|u|_{J_0}\|_2 
  \ge \frac{1}{7\sqrt{\log n}}\|v_0\|_2, 
\end{equation}
where the last inequality holds by Lemma~\ref{C:uj0}.

On the other hand, we can approximate $u$ by $u_0$ as
\begin{equation}                \label{u u0}
  \|u|_{\Lambda}\|_2 
  \le \|u|_{\Lambda} - u_0|_{\Lambda}\|_2 + \|u_0|_{\Lambda}\|_2.
\end{equation}
Since $\L \subset J$ and using Lemma~\ref{L:uj}, we have
$$
\|u|_{\Lambda} - u_0|_{\Lambda}\|_2 \le 2.4\e\|v_0\|_2
$$
Furthermore, by definition \eqref{v0 x0} of $v_0$, we have $v_0|_\Lambda = 0$. 
So, by Part 1 of Proposition~\ref{P:cons}, 
$$
\|u_0|_{\Lambda}\|_2 \le 2.03 \e \|v_0\|_2.
$$
Using the last two inequalities in \eqref{u u0}, we conclude that 
$$
\|u|_{\Lambda}\|_2 \le 4.43 \e \|v_0\|_2.
$$
This is a contradiction to~(\ref{E:ubig}) 
so long as $\e \leq
 0.03 /  \sqrt{\log n}$. 
This proves Theorem~\ref{T:it}.
\qed

\section{Implementation and empirical performance of ROMP}
\label{s: implementation}

\subsection{Running time}

The Identification step of ROMP, i.e. selection of the subset $J$, 
can be done by {\em sorting} the coordinates of $u$ in the nonincreasing order 
and selecting $n$ biggest.
Many sorting algorithms such as Mergesort or Heapsort provide running times 
of $O(d\log d)$. 

The Regularization step of ROMP, i.e. selecting $J_0 \subset J$, 
can be done fast by observing that $J_0$ is an {\em interval} 
in the decreasing rearrangement of coefficients. 
Moreover, the analysis of the algorithm shows that instead of 
searching over all intervals $J_0$, it suffices
to look for $J_0$ among {\em $O(\log n)$ consecutive intervals}
with endpoints where the magnitude of coefficients decreases by a factor of $2$.
(these are the sets $A_k$ in the proof of Lemma~\ref{L:reg}).
Therefore, the Regularization step can be done in time $O(n)$.

In addition to these costs, the $k$-th iteration step of ROMP involves
{\em multiplication} of the $d \times N$ matrix $\Phi^*$ by a vector, 
and solving the {\em least squares problem} with the $N \times |I|$ matrix 
$\Phi_I$, where $|I| \le 2k \le 2n$. 
For unstructured matrices, these tasks can be done in time 
$dN$ and $O(n^2 N)$ respectively. 
Since the submatrix of $\Phi$ when restricted 
to the index set $I$ is near an isometry, using an iterative method 
such as the Conjugate Gradient Method allows us to solve the least 
squares method in a constant number of iterations (up to a specific 
accuracy.) Using such a method then reduces the time of solving the
least squares problem to just $O(nN)$. Thus in the cases where ROMP
terminates after a fixed number of iterations, the total time to solve
all required least squares problems would be just $O(nN)$. 
For structured matrices, such as
partial Fourier, these times can be improved even more using fast multiply
techniques. 

In other cases, however, ROMP may need more than a constant number
of iterations before terminating, say the full $O(n)$ iterations. In this case, 
it may be more efficient to maintain the QR factorization of $\Phi_I$ and use the Modified 
Gram-Schmidt algorithm. With this method, solving all the least squares problems
takes total time just $O(n^2 N).$ However, storing the QR factorization is quite costly, so
in situations where storage is limited it may be best to use the iterative methods mentioned above.

ROMP terminates in at most $2n$ iterations. Therefore, for unstructured 
matrices using the methods mentioned above and in the interesting regime $N \ge \log d$, 
{\em the total running time of ROMP is O(dNn)}. This is the same bound as for OMP
\cite{TG}. 


\subsection{Non-sparse signals}

In many applications, one needs to recover a signal $v$ which is not sparse
but close to being sparse in some way. Such are, for example, compressible
signals, whose coefficients decay at a certain rate (see \cite{Do},
\cite{CT}).
To make ROMP work for such signals, one can replace the stopping criterion
of exact recovery $r=0$ by ``repeat $n$ times or until $r=0$, whichever
occurs first''. Note that we could amend the algorithm for sparse signals
in this way as well, allowing for a specific level of accuracy to be attained
before terminating.

We recently proved that ROMP
is stable and guarantees approximate recovery of non-sparse signals with noisy measurements; this
will be discussed in a forthcoming paper.

\subsection{Experiments}

This section describes our experiments that illustrate the signal recovery power of ROMP.
We experimentally examine how many measurements $N$ are necessary to recover various kinds of $n$-sparse
signals in $\R^d$ using ROMP. 
We also demonstrate that the number of iterations ROMP needs to recover a sparse signal is
in practice at most linear the sparsity. 

First we describe the setup of our experiments. For many values of the ambient dimension $d$, 
the number of measurements $N$, and the sparsity $n$, we reconstruct random signals using ROMP.
For each set of values, we generate an $N \times d$ Gaussian measurement matrix $\Phi$ and then perform $500$ independent trials. The results
we obtained using Bernoulli measurement matrices were very similar.
In a given trial, we generate an $n$-sparse signal $v$ in one of two ways. In either case, we first select the support of the signal by choosing $n$ components uniformly 
at random (independent from the measurement matrix $\Phi$). In the cases where we wish to generate flat signals, we then  set these components to one.\footnote{Our work as well as the analysis of Gilbert and Tropp~\cite{TG}
show that this is a challenging case for ROMP (and OMP). } In the cases where we wish to generate sparse compressible signals, we set the $i^{th}$ component of the support to plus or minus $i^{-1/p}$ for a specified value of $0 < p < 1$. We then execute ROMP with the measurement
vector $x = \Phi v$.  

\begin{figure}[ht] 
  \includegraphics[width=0.8\textwidth,height=3.2in]{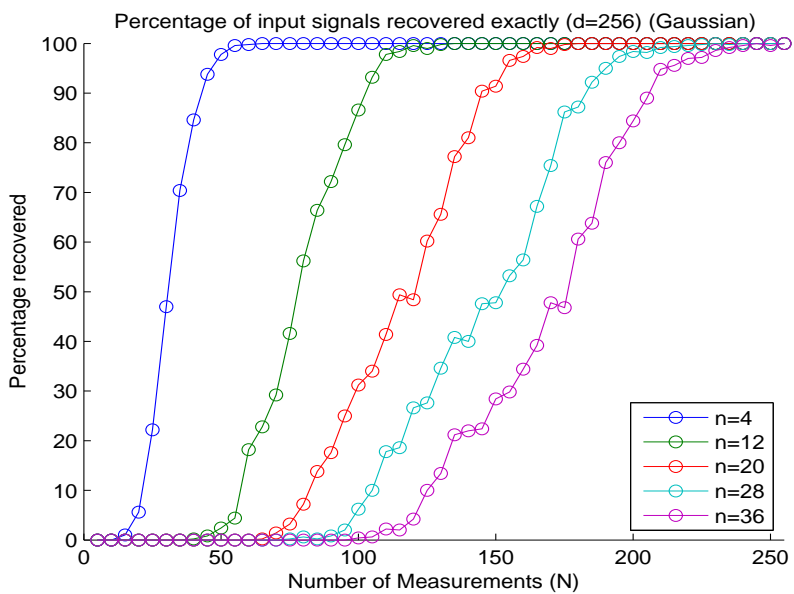}
  \caption{The percentage of sparse flat signals exactly recovered by ROMP as a function of the number of measurements $N$ in dimension $d=256$ for various levels of sparsity $n$.}\label{fig:percent}
\end{figure}

Figure~\ref{fig:percent} depicts the percentage (from the $500$ trials) of sparse flat signals that were 
reconstructed exactly. This plot was generated with $d = 256$ for various levels of sparsity $n$. The horizontal axis represents the number of measurements $N$, and the vertical
axis represents the exact recovery percentage. We also performed this same test for sparse compressible signals and found the results very similar to those in Figure~\ref{fig:percent}.
Our results show that performance of ROMP is very similar to that of OMP which can be found in \cite{TG}. 

\begin{figure}[ht] 
  \includegraphics[width=0.8\textwidth,height=3.2in]{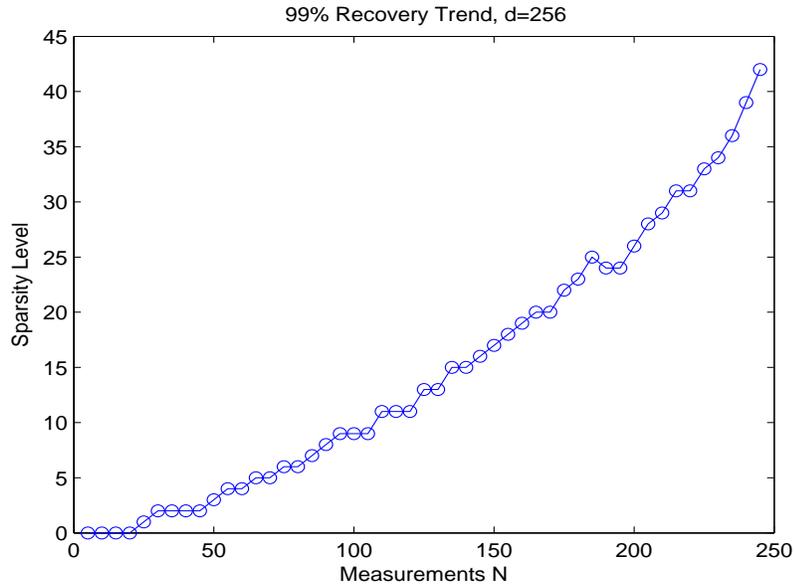}
  \caption{The $99\%$ recovery limit as a function of the sparsity $n$ and the number of measurements $N$ for sparse flat signals.}\label{fig:99}
\end{figure}

Figure~\ref{fig:99} depicts a plot of the values for $N$ and $n$ at which $99\%$ of sparse flat signals are recovered exactly. This plot was generated with $d=256$. The horizontal axis represents the number of measurements $N$, and the vertical axis the sparsity level $n$. 

\begin{figure}[ht] 
  \includegraphics[width=0.8\textwidth,height=3.2in]{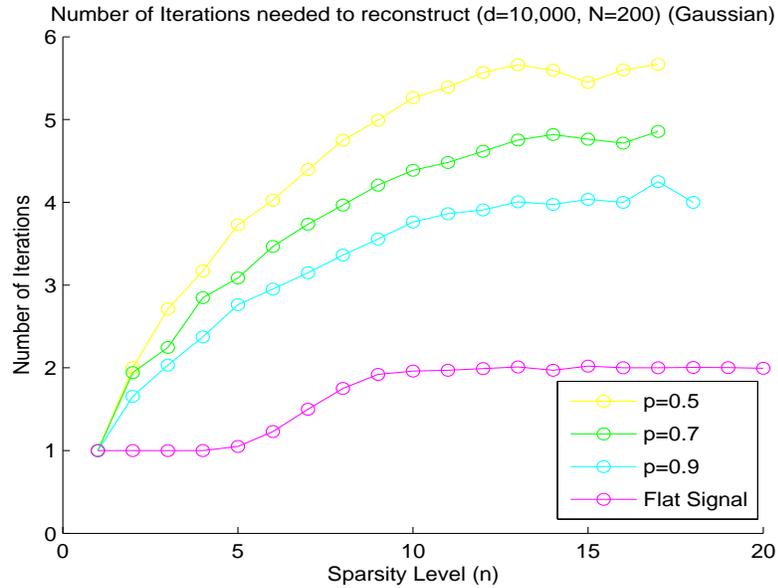}
  \caption{The number of iterations executed by ROMP as a function of the sparsity $n$ in dimension $d=10,000$ with $N=200$.}\label{fig:itsROMP}
\end{figure}

Theorem~\ref{T:main} guarantees that ROMP runs with at most $O(n)$ iterations. Figure~\ref{fig:itsROMP} depicts the number of iterations executed by ROMP for $d=10,000$ and $N=200$. ROMP was 
executed under the same setting as described above for sparse flat signals as well as sparse compressible signals for various values of $p$, and the number of iterations
in each scenario was averaged over the $500$ trials. These averages were plotted against
the sparsity of the signal. As the plot illustrates, only $2$
iterations were needed for flat signals even for sparsity $n$ as high as $40$. The plot also demonstrates that the number of iterations needed for sparse compressible
is higher than the number needed for sparse flat signals, as one would expect. The plot suggests that for smaller
values of $p$ (meaning signals that decay more rapidly) ROMP needs more iterations. However it shows that even in the case of $p=0.5$, only $6$ iterations are needed even for sparsity $n$ as high as $20$.

\end{document}